\documentclass[12pt]{article}
\usepackage{amssymb,latexsym,pstcol,pst-plot,amsthm,amsmath,hyperref}

\newcommand{\cUS}{{\cal US}}

\newcommand{\ben}{\begin{enumerate}}
\newcommand{\een}{\end{enumerate}}
\newcommand{\ble}{\begin{lem}}
\newcommand{\ele}{\end{lem}}
\newcommand{\bth}{\begin{thm}}
\renewcommand{\eth}{\end{thm}}
\newcommand{\bpr}{\begin{prop}}
\newcommand{\epr}{\end{prop}}
\newcommand{\bco}{\begin{cor}}
\newcommand{\eco}{\end{cor}}
\newcommand{\bcon}{\begin{conj}}
\newcommand{\econ}{\end{conj}}
\newcommand{\bde}{\begin{defn}}
\newcommand{\ede}{\end{defn}}
\newcommand{\bex}{\begin{exa}}
\newcommand{\eex}{\end{exa}}
\newcommand{\barr}{\begin{array}}
\newcommand{\earr}{\end{array}}
\newcommand{\btab}{\begin{tabular}}
\newcommand{\etab}{\end{tabular}}
\newcommand{\beq}{\begin{equation}}
\newcommand{\eeq}{\end{equation}}
\newcommand{\bea}{\begin{eqnarray*}}
\newcommand{\eea}{\end{eqnarray*}}
\newcommand{\bce}{\begin{center}}
\newcommand{\ece}{\end{center}}
\newcommand{\bpi}{\begin{picture}}
\newcommand{\epi}{\end{picture}}
\newcommand{\bpp}{\begin{picture}}
\newcommand{\epp}{\end{picture}}
\newcommand{\bfi}{\begin{figure} \begin{center}}
\newcommand{\efi}{\end{center} \end{figure}}
\newcommand{\bprf}{\begin{proof}}
\newcommand{\eprf}{\end{proof}\medskip}
\newcommand{\capt}{\caption}
\newcommand{\bsl}{\begin{slide}{}}
\newcommand{\esl}{\end{slide}}
\newcommand{\bfr}{\begin{frame}}
\newcommand{\efr}{\end{frame}}

\newcommand{\hqed}{\hfill \qed}
\newcommand{\hqedm}{\hfill \qed \medskip}

\newcommand{\hso}[1]{\hspace{-1pt}}

\newcommand{\sbs}{\subset}
\newcommand{\sbe}{\subseteq}

\newcommand{\oph}{\hat{\oplus}}

\newcommand{\ptn}{\vdash}

\def\<{\langle}
\def\>{\rangle}

\newcommand{\spn}[1]{\langle{#1}\rangle}

\newcommand{\si}{\sigma}

\newcommand{\1}{{\bf 1}}

\newcommand{\3}{{\bf 3}}

\newcommand{\bbQ}{{\mathbb Q}}

\newcommand{\cA}{{\cal A}}

\newcommand{\cB}{{\cal B}}

\newcommand{\cE}{{\cal E}}

\newcommand{\cR}{{\cal R}}

\newcommand{\cT}{{\cal T}}

\newcommand{\fS}{{\mathfrak S}}

\newtheorem{thm}{Theorem}[section]
\newtheorem{prop}[thm]{Proposition}
\newtheorem{cor}[thm]{Corollary}
\newtheorem{lem}[thm]{Lemma}
\newtheorem{conj}[thm]{Conjecture}
\newtheorem{exa}[thm]{Example}

\begin{document}
\pagestyle{plain}

\title{Partitions, rooks, and symmetric functions in noncommuting variables
}
\author{Mahir Bilen Can\\[-5pt]
\small Department of Mathematics, Tulane University\\[-5pt]
\small New Orleans, LA 70118, USA, \texttt{mcan@tulane.edu}\\
and\\
Bruce E. Sagan\footnote{Work partially done while a Program Officer at NSF.  The views
    expressed are not necessarily those of the NSF.}\\[-5pt]
\small Department of Mathematics, Michigan State University,\\[-5pt]
\small East Lansing, MI 48824-1027, USA, \texttt{sagan@math.msu.edu}}

\date{\today\\[10pt]
{\it This paper is dedicated to Doron Zeilberger on the occasion of his 60th
  birthday.  His enthusiasm for combinatorics has been an inspiration to us all.}
	\begin{flushleft}
	                                       \\[5pt]
	\end{flushleft}}

\maketitle

\begin{abstract}
Let $\Pi_n$ denote the set of all set partitions of $\{1,2,\ldots,n\}$.  We
consider two subsets of $\Pi_n$, one connected to rook theory and one
associated with symmetric functions in noncommuting variables.  Let
$\cE_n\sbe\Pi_n$ be the subset of all partitions corresponding to an
extendable rook (placement) on the upper-triangular board, $\cT_{n-1}$.  Given
$\pi\in\Pi_m$ and $\si\in\Pi_n$, define their {\it slash product\/} to be
$\pi|\si=\pi\cup(\si+m)\in\Pi_{m+n}$ where $\si+m$ is the partition obtained
by adding $m$ to every element of every block of $\si$.  Call $\tau$ {\it
  atomic\/} if it can not be written as a nontrivial slash product and let
$\cA_n\sbe\Pi_n$ denote the subset of atomic partitions.  Atomic partitions
were first defined by Bergeron, Hohlweg, Rosas, and Zabrocki during their
study of $NCSym$, the symmetric functions in noncommuting variables.  We show
that, despite their very different definitions, $\cE_n=\cA_n$ for all $n\ge0$.
Furthermore, we put an algebra structure on the formal vector space generated
by all rook placements on upper triangular boards which makes it isomorphic to
$NCSym$.  We end with some remarks and an open problem.
\end{abstract}

\section{Extendable rooks and atomic partitions}
\label{era}

For a nonnegative integer $n$, let $[n]=\{1,2,\ldots,n\}$.  Let $\Pi_n$ denote the set 
of all set partitions $\pi$ of $[n]$, i.e., $\pi=\{B_1,B_2,\ldots,B_k\}$ with $\uplus_i B_i = [n]$ (disjoint union).  In this case we will write $\pi\ptn[n]$.  The $B_i$ are called {\it blocks\/}. We will often drop set parentheses and commas and just put slashes  between blocks for readability's sake.  Also, we will always write $\pi$ is {\it standard form\/} which means that
\beq
\label{min}
\min B_1 < \min B_2 <\ldots<\min B_k
\eeq
and the elements in each block are listed in increasing order.  For example 
$\pi=136|2459|78\ptn [9]$.  The {\it trivial partition\/} is the unique element of $\Pi_0$, while all other partitions are {\it nontrivial\/}.  

The purpose of this note is to show that two subsets of $\Pi_n$, one connected with rook theory and the other associated to the Hopf algebra $NCSym$ of symmetric functions in noncommuting variables, are actually equal although they have very different definitions.  After proving this result in the current section, we will devote the next to putting an algebra structure on certain rook placements which is isomorphic to $NCSym$.  The final section contains some comments and open questions.

Let us first introduce the necessary rook theory.  A {\it rook (placement)\/}
is an $n\times n$ matrix, $R$, of $0$'s and $1$'s with at most one 1 in every
row and column.  So a permutation matrix, $P$, is just a rook of full rank.  A
{\it board\/} is $\cB\sbe[n]\times[n]$.  We say that $R$ is a {\it rook on
  $\cB$\/} if $R_{i,j}=1$ implies $(i,j)\in\cB$.  In this case we write, by
abuse of notation, $R\sbe \cB$.  A rook $R\sbe\cB$ is {\it extendable in
  $\cB$\/} if there is a permutation matrix $P$ such that $P_{i,j}=R_{i,j}$
for $(i,j)\in\cB$.  For example, consider the upper-triangular board
$\cT_n=\{(i,j)\ :\ i\le j\}$.  The $R\sbe\cT_2$ are displayed in
Figure~\ref{T2}.  Only the third and fifth rooks in Figure~\ref{T2} are
extendable, corresponding to the transposition and identity permutation
matrices, respectively.  Extendability is an important concept in rook theory
because of its relation to the much-studied hit numbers of a
board~\cite[page 163 and ff.\/]{rio:ica}. 

\bfi
$$
\barr{rccccc}
R:
&
\left(
\barr{cc}
0&0\\
0&0
\earr
\right)
&
\left(
\barr{cc}
1&0\\
0&0
\earr
\right)
&
\left(
\barr{cc}
0&1\\
0&0
\earr
\right)
&
\left(
\barr{cc}
0&0\\
0&1
\earr
\right)
&
\left(
\barr{cc}
1&0\\
0&1
\earr
\right)
\\[20pt]
\pi_R:
&
1|2|3
&
12|3
&
13|2
&
1|23
&
123
\earr
$$
\capt{The rooks on $\cT_2$ and their associated partitions}
\label{T2}
\efi

There is a well-known bijection between $\pi\in\Pi_n$ and the rooks
$R\sbe\cT_{n-1}$~\cite[page 75]{sta:ec1}.  Given $R$, define a partition $\pi_R$ by
putting $i$ and $j$ in the same block of $\pi_R$ whenever $R_{i,j-1}=1$.  For
each $R\sbe\cT_2$, the corresponding $\pi_R\in\Pi_3$ is shown in
Figure~\ref{T2}.  Conversely, given $\pi$ we define a rook $R_\pi$ by letting
$(R_\pi)_{i,j}=1$ exactly when $i$ and $j+1$ are adjacent elements in a block
of $\pi$ in standard form.  It is easy to see that the maps $R\mapsto\pi_R$
and $\pi\mapsto R_\pi$ are inverses.  If a matrix has a certain property  then
we will also say that the corresponding partition does, and vice-versa.  Our
first subset of $\Pi_n$ will  be the {\it extendable partitions\/} denoted by 
$$
\cE_n=\{\pi\in\Pi_n\ :\ \mbox{$R_\pi$ is extendable in $\cT_{n-1}$}\}.
$$
So, from Figure~\ref{T2}, $\cE_2=\{13|2,123\}$.

To define our second subset of $\Pi_n$, it is convenient to introduce an operation on partitions.  For a set of integers $B=\{b_1,\ldots,b_j\}$ we let $B+m=\{b_1+m,\ldots,b_j+m\}$.  Similarly, for a partition $\pi=\{B_1,\ldots,B_k\}$ we use the notation $\pi+m=\{B_1+m,\ldots,B_k+m\}$.  If $\pi\in\Pi_m$ and $\si\in\Pi_n$ then define their {\it slash product\/} to be the partition in $\Pi_{m+n}$ given by
$$
\pi|\si=\pi\cup(\si+m).
$$
Call a partition {\it atomic\/} if it can not be written as a slash product of two nontrivial partitions and let
$$
\cA_n=\{\pi\in\Pi_n\ :\ \mbox{$\pi$ is atomic}\}.
$$
Atomic partitions were defined by Bergeron, Hohlweg, Rosas, and Zabrocki~\cite{bhrz:gbp} because of their connection with symmetric functions in noncommuting variables.  We will have more to say about this in Section~\ref{arp}.  

Since $\cE_n$ is defined in terms of rook placements, it will be convenient to have a rook interpretation of $\cA_n$.  Given any two matrices $R$ and $S$, defined their {\it extended direct sum\/} to be
$$
R\oph S = R\oplus (0) \oplus S
$$ 
where $\oplus$ is ordinary matrix direct sum and $(0)$ is the $1\times 1$ zero matrix.  To illustrate,
$$
\left(
\barr{ccc}
a&b&c\\
d&e&f
\earr
\right)
\oph
\left(
\barr{cc}
w&x\\
y&z
\earr
\right)
=
\left(
\barr{cccccc}
a&b&c&0&0&0\\
d&e&f&0&0&0\\
0&0&0&0&0&0\\
0&0&0&0&w&x\\
0&0&0&0&y&z
\earr
\right).
$$
It is clear from the definitions that $\tau=\pi|\si$ if and only if 
$R_\tau=R_\pi\oph R_\si$.  We now have everything we need to prove our first result.

\bth
\label{E=A}
For all $n\ge0$ we have $\cE_n=\cA_n$.
\eth 
\proof
Suppose we have $\tau\in\cE_n$.  Assume, towards a contradiction, that $\tau$ is not atomic so that $\tau=\pi|\si$.  On the matrix level we have $R_\tau=R_\pi\oph R_\si$ where $R_\pi$ is $m\times m$ for some $m$.  We are given that $\tau$ is extendable, so let $P$ be a permutation matrix extending $R_\tau$.  Since $P$ and $R_\tau$ agree above and including the diagonal, the first $m+1$ rows of $P$ must be zero from column $m+1$ on.  But $P$ is a permutation matrix and so each of these $m+1$ rows must have a one in a different column, contradicting the fact that only $m$ columns are available.

Now assume $\tau\in\cA_n$.  We will construct an extension $P$ of $R_\tau$.  Let $i_1,\ldots,i_r$ be the indices of the zero rows of $R_\tau$ and similarly for $j_1,\ldots,j_r$ and the columns.  If $i_k>j_k$ for all $k\in[r]$, then we can construct $P$ by supplementing $R_\tau$ with ones in positions $(i_1,j_1),\ldots,(i_r,j_r)$.

So suppose, towards a contradiction, that there is some $k$ with $i_k\le j_k$.  Now $R_\tau$ must contain $j_k-k$ ones in the columns to the left of column $j_k$.  If $i_k<j_k$, then there are fewer than $j_k-k$ rows which could contain these ones since $R_\tau$ is upper triangular. This is a contradiction.  If $i_k=j_k$, then the $j_k-k$ ones in the columns left of $j_k$ must lie in the first $i_k-k=j_k-k$ rows.  Furthermore, these ones together with the zero rows force the columns to the right of $j_k$ to be zero up to and including row $i_k=j_k$.  It follows that $R_\tau=R_\pi\oph R_\si$ for some $\pi,\si$ with $R_\pi$ being $(i_k-1)\times(i_k-1)$.  This contradicts the fact that $\tau$ is atomic.\hqedm  

Having two descriptions of this set may make it easy to prove assertions about it from one definition which would be difficult to demonstrate if the other were used.  Here is an example.

\bco
Let $R\sbe\cT_n$.  If $R_{1,n}=1$ then $R$ is extendable in $\cT_n$.
\eco
\proof
If $R_{1,n}=1$ then we can not have $R=R_\si\oph R_\tau$ for nontrivial $\si,\tau$.  So $R$ is atomic and, by the previous theorem, $R$ is extendable.\hqedm

\section{An algebra on rook placements and $NCSym$}
\label{arp}

The algebra of symmetric functions in noncommuting variables, $NCSym$, was first studied by Wolf~\cite{wol:sfn} who proved a version of the Fundamental Theorem of Symmetric Functions in this context.  The algebra was rediscovered by Gebhard and Sagan~\cite{gs:csf} who used it as a tool to make progress on Stanley's ($\3+\1$)-free Conjecture for chromatic symmetric functions~\cite{sta:sfg}.  Rosas and Sagan~\cite{rs:sfn} were the first to make a systematic study of the vector space properties of $NCSym$.  Bergeron, Reutenauer, Rosas, and Zabrocki~\cite{brrz:ics} introduced a Hopf algebra structure on $NCSym$ and described its invariants and covariants.  

Let $X=\{x_1,x_2,\ldots\}$ be a countably infinite set of variables which do not commute.  Consider the corresponding ring of formal power series over the rationals $\bbQ\spn{\spn{X}}$.  Let $\fS_m$ be the symmetric group on $[m]$.  Then any $g\in\fS_n$ acts on a monomial $x=x_{i_1} x_{i_2}\cdots x_{i_n}$ by
$$
g(x)=x_{g^{-1}(i_1)} x_{g^{-1}(i_2)}\cdots x_{g^{-1}(i_n)}
$$
where $g(i)=i$ for $i>m$.  Extend this action linearly to $\bbQ\spn{\spn{X}}$.  The {\it symmetric functions in noncommuting variables\/}, $NCsym\sbs\bbQ\spn{\spn{X}}$, are all power series which are of bounded degree and invariant under the action of $\fS_m$ for all $m\ge0$.

The vector space bases of $NCSym$ are indexed by set partitions.  We will be particularly interested in a basis which is the analogue of the power sum basis for ordinary symmetric functions.  Given a monomial $x=x_{i_1} x_{i_2}\cdots x_{i_n}$, there is an associated set partition $\pi_x$ where $j$ and $k$ are in the same block of $\pi_x$ if and only if $i_j=i_k$ in $x$, i.e., the indices in the $j$th and $k$th positions are the same.  For example, if $x=x_3 x_5 x_2 x_3 x_3 x_2$ then $\pi_x=145|2|36$.  The {\it power sum symmetric functions in noncommuting variables\/} are defined by
$$
p_\pi=\sum_{x\ :\ \pi_x\ge\pi} x,
$$
where $\pi_x\ge\pi$ is the partial order in the lattice of partitions, so $\pi_x$ is obtained by merging blocks of $\pi$.  Equivalently, $p_\pi$ is the sum of all monomials where the indices in the $j$th and $k$th places are equal if $j$ and $k$ are in the same block of $\pi$, but there may be other equalities as well.  To illustrate,
$$
p_{13|2}=x_1 x_2 x_1 + x_2 x_1 x_2 + \cdots + x_1^3 + x_2^3 + \cdots.
$$

Note that, directly from the definitions,
\beq
\label{p}
p_{\pi|\si} = p_\pi p_\si.
\eeq
Using this property, Bergeron, Hohlweg, Rosas, and Zabrocki~\cite{bhrz:gbp} proved the following result which will be useful for our purposes.
\bpr[\cite{bhrz:gbp}]
\label{p_pi}
As an algebra, $NCSym$ is freely generated by the $p_\pi$ with $\pi$ atomic.\hqed
\epr

Let 
$$
\cR=\{R\sbe\cT_n\ :\ n\ge-1\},
$$
where there is a single rook on $\cT_{-1}$ called the {\it unit rook\/} and denoted $R=1$ (not to be confused with the empty rook on $\cT_0$).  We extend the bijection between set partitions and rooks on upper triangular boards by letting the unit rook correspond to the empty partition.  Consider the vector space $\bbQ\cR$ of all formal linear combinations of rooks in $\cR$.  By both extending $\oph$ linearly and letting the unit rook act as an identity, the operation of extended direct sum can be considered as a product on this space.  It is easy to verify that this turns $\bbQ\cR$ into an algebra.
\bpr
\label{R_pi}
As an algebra, $\bbQ\cR$ is freely generated by the $R_\pi$ with $\pi$ atomic.
\epr
\proof
A simple induction on $n$ shows  that any $\tau\in\Pi_n$ can be uniquely factored as $\tau=\pi_1|\pi_2|\cdots|\pi_t$ with the $\pi_i$ atomic.  From the remark just before Theorem~\ref{E=A}, it follows that each $R_\tau$ can be uniquely written as a product of  atomic $R_\pi$'s.  Since the set of all $R_\tau$ forms a vector space basis, the atomic $R_\pi$ form a free generating set.\hqedm

Comparing Propositions~\ref{p_pi} and~\ref{R_pi} as well as the remark before Theorem~\ref{E=A} and equation~\ref{p}, we immediately get the desired isomorphism.
\bth
The map $p_\pi\mapsto R_\pi$ is an algebra isomorphism of $NCSym$ with $\bbQ\cR$.\hqed
\eth

\section{Remarks and an open question}

\subsection{Unsplittable partitions}

Bergeron, Reutenauer, Rosas, and Zabrocki~\cite{brrz:ics} considered another free generating set for $NCSym$ which we will now describe.  A {\it restricted growth function of length $n$\/} is a sequence of positive integers $r=a_1a_2\ldots a_n$ such that
\ben
\item $a_1=1$, and
\item $a_i\leq 1+\max\{a_1,\ldots,a_{i-1}\}$ for $2\le i\le n$.
\een
Let $RG_n$ denoted the set of restricted growth functions of length $n$.
There is a well-known bijection between $\Pi_n$ and $RG_n$~\cite[page 34]{sta:ec1} as
follows.  Given $\pi\in\Pi_n$ we define $r_\pi$ by $a_i=j$ if and only if
$i\in B_j$ in $\pi$.  For example, if $\pi=124|36|5$ then $r_\pi=112132$.  It
is easy to see that having $\pi$ in standard form makes the map well defined.
And the reader should have no trouble constructing the inverse. 

Define the {\it split product\/} of $\pi\in\Pi_m$ and $\si\in\Pi_n$ to be $\tau=\pi\circ\si\in\Pi_{m+n}$ where $\tau$ is the uniqe partition such that $r_\tau=r_\pi r_\si$ (concatenation).  To illustrate, if $\pi$ is as in the previous paragraph and $\si=13|2$ then $r_\pi r_\si = 112132121$ and so $\pi\circ\si=12479|368|5$.  This is not Bergeron et al.'s original definition, but it is equivalent.  Now define $\tau$ to be {\it unsplitable\/} if it can not be written as a split product of two nontrivial partitions.  (Bergeron et al.\ used the term "nonsplitable" which is not a typical English word.)  Let $\cUS_n\sbe\Pi_n$ be the subset of unsplitable partitions.  So $\cUS_2=\{1|2|3, 1|23\}$.

Perhaps the simplest basis for $NCSym$ is the one gotten by symmetrizing a monomial.  Define the {\it monomial symmetric functions in noncommuting variables\/} to be
$$
m_\pi=\sum_{x\ :\ \pi_x=\pi} x.
$$
So now indices in a term of $m_\pi$ are equal precisely when their positions are in the same block of $\pi$.  For example,
$$
m_{13|2}=x_1 x_2 x_1 + x_2 x_1 x_2 + \cdots.
$$
The following is a more explicit version of Wolf's original result~\cite{wol:sfn}.
\bpr[\cite{brrz:ics}]
\label{m_pi}
As an algebra, $NCSym$ is freely generated by the $m_\pi$ with $\pi$ unsplitable.\hqed
\epr

Comparing Propositions~\ref{p_pi} and~\ref{m_pi} we see that $|\cA_n|=|\cUS_n|$ for all $n\ge0$ where $|\cdot|$ denotes cardinality.  (Although they are not the same set as can be seen by our computations when $n=2$.)  It would be interesting to find a bijective proof of this result.

\subsection{Hopf structure}

Thiem~\cite{thi:brr} found a connection between $NCSym$ and unipotent
upper-triangular zero-one matrices using supercharacter theory.   This work
has very recently been extended using matrices over any field and a colored
version of $NCSym$ during a workshop at the American Institute of
Mathematics~\cite{ber:pc}. 
This approach gives an 
isomorphism even at the Hopf algebra level.

\end{document}